\documentclass[preprint,3p,12pt,pdf]{elsarticle}

\usepackage{lineno,hyperref}
\modulolinenumbers[5]

\usepackage{hyperref}

\usepackage{epstopdf}

\usepackage{tikz}
\usetikzlibrary{arrows}

\usepackage{pst-node}
\usepackage{tikz-cd}

\usepackage[shellescape]{gmp}

\usepackage{mathrsfs}
\usepackage{amssymb}
\usepackage{amsfonts}
\usepackage{latexsym}
\usepackage{mathtools}
\usepackage{xcolor}
\usepackage{wrapfig}
\usepackage{floatflt}

\usepackage{mathtools}
\usepackage{extarrows}

\usepackage{graphicx}

\usepackage{subcaption}

\def\lb{\label}

\newcommand{\er}[1]{\textrm{(\ref{#1})}}

\def\a{\alpha}

\def\G{\Gamma}

\def\m{\mu}            
            
\def\r{\rho}

\def\O{\Omega}

\def\ve{\varepsilon}       \def\vp{\varphi}

    \def\N{{\mathbb N}}   

\def\lt{\biggl}                  \def\rt{\biggr}
               \def\wt{\widetilde}

\let\ge\geqslant                 \let\le\leqslant

\def\iy{\infty}

\def\el2{\ell^{\,2}}             \def\1{1\!\!1}

\def\Re{\mathop{\mathrm{Re}}\nolimits}

\let\ge\geqslant
\let\le\leqslant

\newcommand{\ca}{\begin{cases}}
	\newcommand{\ac}{\end{cases}}
\newcommand{\ma}{\begin{pmatrix}}
	\newcommand{\am}{\end{pmatrix}}
\renewcommand{\[}{\begin{equation}}
	\renewcommand{\]}{\end{equation}}
\def\eq{\begin{equation}}
	\def\qe{\end{equation}}

\begin{document}
	
	\begin{frontmatter}

		\title{Complete left tail asymptotic for branching processes in random environments}

		\date{\today}

		\author
		{Anton A. Kutsenko}
	
\address{University of Hamburg, MIN Faculty, Department of Mathematics, 20146 Hamburg, Germany; email: akucenko@gmail.com}
	
	\begin{abstract}
		Recently, the complete left tail asymptotic for the density of the {\it martingale limit} of the classical Galton-Watson process has been derived. The derivation is based on the properties of a special function (whose inverse Fourier transform is the density) satisfying a Poincar\'e-type functional equation. However, the corresponding special function does not satisfy the Poincar\'e-type functional equation for the processes in random environments. This makes it impossible to apply the developed asymptotic methods. In this instance, even computing the special function and the density becomes extremely challenging. We propose a method that extracts the terms of the asymptotic series from another limit distribution {\it ratios of probabilities of rare events}. The generating functions for this limit distribution satisfy the generalized Schr\"oder-type functional equations that greatly simplify the calculations. As in the classical case of constant environments, the obtained series converges everywhere. Possible interesting connections with exotic fractal measures that Robert Strichartz considers are also discussed.
	\end{abstract}


	
\end{frontmatter}


{\section{Introduction}\lb{sec0}}

About the classical branching processes with constant probability-generating function, one may read \cite{H}, where some limit theorems and asymptotics are discussed. For branching processes in a random environment, the probability-generating function is chosen at each step randomly from some family of functions. The study of such processes dates back to \cite{AK71} and continues actively to this day. For example, in \cite{BHS}, the authors consider a random environment oscillating between `healthy' and `harsh' states. Due to such processes' complexity, only the dominant behavior is usually estimated. We will try to find a complete description of limit distributions.    

If the mathematical expectation of the reproduction of offspring is the same for each step, then the population grows exponentially with the exponent equal to the expectation. Thus, one can define a {\it martingale limit}, which shows the spread of population growth compared to the basic exponential growth. The density of the {\it martingale limit} can be computed as a Fourier transform of some special function. The special function and its Fourier transform are quite complex, so any asymptotic of the density is welcome. We will discuss only the so-called Schr\"oder case of branching processes, where an individual's survival probability at each time step is non-zero. Let us recall some results on the classical case of constant environment. For the small argument of the density, the first asymptotic term was analysed in \cite{D1}, where precise estimates were obtained. Then in \cite{BB1}, the authors found an explicit form of the first asymptotic term. After that, in \cite{K24}, the complete left tail asymptotic was derived. All the asymptotic terms were expressed explicitly through the so-called one-periodic Karlin-McGregor function introduced in \cite{KM1} and \cite{KM2}. Moreover, it was shown that the series converges everywhere, not only for small arguments of the density. In random environments, defining a single one-periodic Karlin-McGregor function is impossible. Other approaches must be developed, which will be discussed below.

Let $X_t$ be a Galton-Watson branching process in random environments
\[\lb{000}
 X_{t+1}=\sum_{i=1}^{X_t}\xi_{it},\ \ \ X_0=1,\ \ \ t\in\N\cup\{0\},
\]
where all $\xi_{it}$  are independent and, for any fixed $t$, identically-distributed natural number-valued random variables with offspring probability-generating function 
\[\lb{001}
 P_r(z)=p_{1r}z+p_{2r}z^2+p_{3r}z^3+....
\]
Here $p_{jr}$ are the probabilities of producing $1$, $2$, ... offspring by individuals in a generation - they are all non-negative and their sum is $\sum_{j=1}^{+\iy} p_{jr}=1$. Each time $t$, the parameter $r$ is chosen randomly from some probabilistic space $(R,\O,\m)$. The function $P_r(z)$ is the same for all individuals at a fixed time $t$. At the next time $t+1$, the function $P_{\wt r}(z)$, again the same for all individuals at this time $t+1$, can differ from $P_r(z)$ because $\wt r$ is randomly chosen. This is a typical situation when the external environment randomly changes in some ranges of possible scenarios $(R,\O,\m)$, affecting the reproduction function of offspring $P_r(z)$. 

We assume that the expectation 
\[\lb{002}
 E:=P_r'(1)=p_{1r}+2p_{2r}+3p_{3r}+...<+\iy
\]
is finite and does not depend on $r$. Then the family grows exponentially as $E^t$ and one can define the {\it martingale limit} $\lim_{t\to+\iy}E^{-t}X_t$. In the classical case of the constant environment, when we have one probability-generating function $P(z)$, the density of the {\it martingale limit} can be computed as
\[\lb{003}
 p(x)=\frac1{2\pi}\int_{-\iy}^{+\iy}\Pi(\mathbf{i}y)e^{\mathbf{i}yx}dy,
\ \ \ {\rm where}\ \ \ 
 \Pi(z)=\lim_{t\to+\iy}\underbrace{P\circ...\circ P}_{t}(1-\frac{z}{E^{t}}),
\]
see, e.g., \cite{D1}. By definition, it is seen that $\Pi$ satisfies the so-called Poincar\'e-type functional equation
\[\lb{004}
 P(\Pi(z))=\Pi(Ez),\ \ \ \Pi(0)=1,\ \ \ \Pi'(0)=-1.
\]
This fact allows one to define the so-called one-periodic Karlin-McGregor function, and use it to derive the complete left-tail asymptotic series that converges everywhere, not only for small arguments, see \cite{K24}. This result completes the research that began with excellent work \cite{BB1}, where the first term of asymptotics was obtained.

For the branching process in random environments, instead of \er{003}, one should consider all the possible compositions of the probability-generating functions. Thus,
\[\lb{100}
 \Pi(z)=\lim_{t\to+\iy}\int_R...\int_RP_{r_1}\circ...\circ P_{r_t}\lt(1-\frac{z}{E^t}\rt)\m(dr_t)...\m(dr_1),
\]
which then determines the density in the standard way
\[\lb{101}
 p(x)=\frac1{2\pi}\int_{-\iy}^{+\iy}\Pi(\mathbf{i}y)e^{\mathbf{i}yx}dy=\frac1{2\pi}\int_{-\iy}^{+\iy}\Pi(\ve+\mathbf{i}y)e^{(\ve+\mathbf{i}y)x}dy,
\]
for any $\ve>0$. The corresponding derivations repeat the steps from \cite{D1} for the constant environments in the way explained in the Appendix of \cite{K24a}, made for random environments. 

In practice, other models of time-varying branching processes, e.g., Crump-Mode-Jagers ones, may lead to integral equations different from \er{100}. However, some connections between models are probably possible if we find a common interpretation of time scaling in {\it martingale limit} and infected individuals' life periods - the qualitative form of profiles in Figs. \ref{fig1} and Fig. 2 in \cite{PMPWBMMB} are similar despite the different meanings.

{\bf Remark 1.} For the computation of $\Pi(z)$ in \er{003}, there exist fast procedures developed in, e.g., \cite{Du1} and \cite{K}, providing exponentially fast convergence. However, even if $R$ consists only of two points, say $\{0\}$ and $\{1\}$, with different probability-generating functions $P_0(z)$ and $P_1(z)$, the computation of \er{100} can be difficult. Going into details, let us define
\[\lb{101aa}
 \Pi_{\r}(z)=\lim_{t\to+\iy}P_{r_1}\circ...\circ P_{r_t}\lt(1-\frac{z}{E^t}\rt),\ \ \ \r=\frac{r_1}2+\frac{r_2}{2^2}+\frac{r_3}{2^3}+....
\]
Then
\[\lb{101ab}
 \Pi(z)=\int_0^1\Pi_{\r}(z)\pmb{\m}(d\r),
\]
where the measure $\pmb{\m}$ is the defined on diadic intervals as
$$
  \pmb{\m}([\r,\r+2^{-t}))=\prod_{j=1}^t\m(\{r_j\}),\ \ \ \r=\frac{r_1}2+\frac{r_2}{2^2}+\frac{r_3}{2^3}+...+\frac{r_t}{2^t}.
$$
If $\m(\{0\})\ne\m(\{1\})$ then $\pmb{\m}$ is a fractal, singular measure. Robert S. Strichartz in \cite{S}  considered some complex properties of such measures. Subsequent results, such as the Fourier transform of such measures and the explicit calculation of integrals, were obtained in \cite{K23a}. If  $\m(\{0\})=\m(\{1\})$ then $\pmb{\m}$ is the standard Lebesgue measure. Nevertheless, the computation of \er{101ab} is still hard, since $\Pi_{\r}(z)$ depends on $\r$ irregularly, with no continuity. Even if the computation of a single $\Pi_{\r}(z)$ is not very difficult for some cases, the parameter $\r$ takes a continuum number of values. 

It is easy to see that there are no analogs of the first formula in \er{004} for $\Pi(z)$ defined in \er{100}. This does not allow us to introduce the one-periodic Karlin-McGregor function and use the methods from \cite{K24} to obtain the left tail asymptotics. We will go another way.

The {\it ratios of probabilities of rare events} or {\it relative limit densities of the number of descendants} are defined by
\[\lb{102}
 \vp_{nj}=\lim_{t\to+\iy}\frac{\mathbb{P}(X_t=n)}{\mathbb{P}(X_t=j)},\ \ {\rm under\ condition}\ \ \mathbb{P}(X_0=j)=1.
\]
The corresponding generating functions
\[\lb{103}
 \Phi_j(z)=\sum_{n=j}^{+\iy}\vp_{nj}z^n
\]
satisfy the generalized Schr\"oder type functional equation
\[\lb{104}
 \int_R\Phi_j(P_r(z))\m(dr)=q_j\Phi_j(z),\ \ \ q_j:=\int_Rp_{r1}^j\m(dr),\ \ \ \Phi_j(z)\sim z^j,\ z\to0.
\]
The derivation of \er{104} repeats the arguments from \cite{K24a}, see Appendix in that paper. In \cite{K24a}, I promised a potential interest of unusual {\it ratios of probabilities of rare events} in the analysis of the more standard {\it martingale limit}. However, at that moment, the connection between these two limit distributions was unclear, especially in the case of random environments. Very soon, everything falls into place. 

Substituting \er{103} into \er{104}, we obtain the recurrence formula for the computation $\vp_{n,j}$, namely
\begin{multline}\lb{105}
 \vp_{1j}=...=\vp_{j-1,j}=0,\ \ \ \ \ \vp_{jj}=1,\ \ \ \ \ \vp_{nj}=\frac{\sum_{m=1}^{n-1}q_{nm}\vp_{mj}}{q_{jj}-q_{nn}},\ n>j,\\ q_{nm}=\sum_{i_1+...+i_m=n}\int_Rp_{ri_1}...p_{ri_{m}}\m(dr),
\end{multline}
where all $i_k\ge1$. In particular,
\[\lb{105a}
 q_{nn}=q_n=\int_Rp_{r1}^n\m(dr),\ \ \ q_{n,n-1}=(n-1)\int_Rp_{r1}^{n-2}p_{r2}\m(dr).
\]
Suppose the probability space $R$ consists of a finite number of points and the corresponding polynomials $P_r(z)$ have small degrees. Then the computation of $\vp_{nj}$ seems not very difficult, at least for moderate $n$ and $j$. It is crucial that $\Pi(z)$ in \er{100} does not satisfy any reasonable functional equation, but functions $\Phi_j(z)$ in \er{103} satisfy very nice functional equations \er{104}. But what is the connection between functions $\Pi(z)$ or $p(x)$, see \er{101}, and $\Phi_j(z)$? First, let us define the ``pseudo-inverse function'' related to $\Phi_j(z)$. Namely, consider the formal identity
\[\lb{106}
 z=\sum_{j=1}^{+\iy}b_j\Phi_j(z).
\]
Due to \er{103} and first identities in \er{105}, the coefficients $b_j$ satisfy
\[\lb{107}
 \ma 1 & 0 & 0 & 0 & ... \\
       \vp_{21} & 1 & 0 & 0 & ... \\
       \vp_{31} & \vp_{32} & 1 & 0 & ... \\
       \vp_{41} & \vp_{42} & \vp_{43} & 1 & ... \\
       ... & ... & ... & ... & ...
 \am
 \ma
    b_1 \\ b_2 \\ b_3 \\ b_4 \\ ...
 \am
=
\ma 
  1 \\ 0 \\ 0 \\ 0 \\ ...
\am,
\]
which leads to the recurrence formula
\[\lb{108}
 b_1=1,\ \ \ b_n=-\sum_{m=1}^{n-1}\vp_{nm}b_m.
\]
Another formula for the computation of $b_n$ is based on applying Cramer's rule to \er{107}. It is
\[\lb{109}
 b_1=1,\ \ \ \ \  b_n=(-1)^n\begin{vmatrix}
\vp_{21} & 1 & 0 & 0 & ... \\ 
\vp_{31} & \vp_{32} & 1 & 0 & ... \\ 
\vp_{41} & \vp_{42} & \vp_{43} & 1 & ... \\
... & ... & ... & ... & ... \\
\vp_{n1} & \vp_{n2} & \vp_{n3} & ... & \vp_{n,n-1}
\end{vmatrix},\ \ \ n>1.
\]
Of course, \er{108} and \er{109} give the same result.

{\bf Remark 2.} For the case of constant environments, when we have only one probability-generating function $P(z)$, the functions $\Phi_j(z)$ are determined through the first function explicitly $\Phi_j(z)=\Phi_1(z)^j$. Then the pseudo-inverse identity \er{106} becomes
$$
 z=\Phi_1^{-1}(\Phi_1(z))=\sum_{j=1}^{+\iy}b_j\Phi_j(z),
$$
where $b_j$ are Taylor coefficients of $\Phi_1^{-1}(z)$ at $z=0$. There are no problems with the existence of $\Phi_1^{-1}(z)$ and the convergence of the corresponding Taylor series in some disk, since $\Phi_1(z)\sim z$, $z\to0$ is analytic for $|z|<1$, see, e.g., \cite{K24a}. However, in the case of random environments, the proof of convergence of \er{106} can be cumbersome. 

Let us formulate the main result. For simplicity, let us assume that $R$ consists of a finite number of points and all $P_r(z)$ are (non-trivial) polynomials. This is a technical requirement that may be significantly relaxed. Suppose we have only a finite number of polynomials $P_r(z)$. In that case, we assume that their common {\it critical angle} near $1$ of their Julia sets is strictly greater than $\pi$, see details in \cite{K24a}. We need this for the exponentially fast convergence of some series appearing in the derivation below. Even without understanding the details about the {\it critical angle}, the derivation in the next Section is quite transparent. Another natural assumption is that \er{106} converges in some disk. Then,
\[\lb{vphi_main}
 \vp_{mj}\sim m^{-1-\log_Eq_j}A_j(-\log_Em),\ \ \ m\to+\iy
\]
with some one-periodic functions $A_j(x)$. Extracting leading terms $A_j(x)$ from the values $\vp_{mj}$ computed by \er{105}, one can use them for the computation of the density by
\[\lb{p_main}
 p(x)=\sum_{j=1}^{+\iy}b_jx^{-1-\log_E q_j}A_j(-\log_Ex),
\]
where $b_j$ are given by \er{108} or \er{109}. In \er{122}, we derived much more than provided in \er{vphi_main}. Namely, we derived a complete asymptotic series for $\vp_{mj}$. In particular, the second asymptotic term is $m$ times smaller than the first term. This remark allows us to control the discrepancy in \er{vphi_main}. Moreover, the oscillations $A_j(x)$ are usually small - $A_j(x)$ are almost constants. Phenomena of this type are known in physics and biology when applied to other similar objects, see, e.g., \cite{DIL}, \cite{CG}, \cite{DMZ}. It is called extremely small critical amplitudes. The explanation of this phenomenon is based on the fact that, roughly speaking, the width of the strip, where one-periodic analytic function $A_j(z)$ is defined, is quite large. This is because the {\it critical angle}, discussed above \er{vphi_main}, is larger than $\pi$. Some details are available in, e.g., \cite{K}.

Summarizing the abovementioned observations, one can suggest a good approximation of $p(x)$, see \er{p_main}. We compute 
\[\lb{approx1}
 A_{jM}:=M^{1+\log_Eq_j}\vp_{Mj}
\]
for some large $M$ and then take
\[\lb{approx2}
 p_{JM}(x)=\sum_{j=1}^{J}b_jA_{jM}x^{-1-\log_E q_j}
\]
for some large $J$. Usually, $J$ is not too large, since $A_{jM}$ decays to $0$ very fast in $j$.

In the next Section, we derive \er{vphi_main} and \er{p_main}. Then, we consider examples and conclude.

{\section{Derivation of \er{vphi_main} and \er{p_main}}\lb{sec1}}

\subsection{Analogs of multiplicatively periodic Karlin-McGregor functions.} Using  \er{106}, we obtain
\[\lb{110a}
 P_{r_1}\circ...\circ P_{r_t}\lt(1-\frac{z}{E^t}\rt)=\sum_{j=1}^{+\iy}b_j\Phi_j\lt(P_{r_1}\circ...\circ P_{r_t}\lt(1-\frac{z}{E^t}\rt)\rt),
\] 
which, after taking the integrals and the limit, see \er{100}, turn into
\[\lb{110}
 \Pi(z)=\sum_{j=1}^{+\iy}b_jK_j(z),
\]
where
\[\lb{111}
 K_j(z)=\lim_{t\to+\iy}\int_R...\int_R\Phi_j\lt(P_{r_t}\circ...\circ P_{r_1}\lt(1-\frac{z}{E^t}\rt)\rt)\m(dr_1)...\m(dr_t).
\]
Due to \er{104}, the functions $K_j(z)$ are multiplicatively periodic
\[\lb{112}
 K_j(Ez)=q_jK_j(z).
\]
Identity \er{112} allows us to determine the one-periodic function
\[\lb{113}
 L_j(z):=q_j^{-z}K_j(E^z),\ \ \ L_j(z+1)=L_j(z),
\]
which can be expanded into the Fourier series
\[\lb{114}
 L_j(z)=\sum_{n=-\iy}^{+\iy}L_{nj}e^{2\pi\mathbf{i}nz},
\] 
with some coefficients $L_{nj}$. By analogy with the results from \cite{K24a}, the coefficients $L_{nj}$ decay in $n$ exponentially fast. The decaying rate depends on the {\it critical angle} mentioned in the Introduction Section. Such a fast decay allows one to change the order of summation of the series and the order of summation and integration, which we will use below. Using \er{113}, we obtain
\[\lb{115}
 K_j(z)=z^{\log_Eq_j}L_j(\log_Ez),
\]
which with \er{114} leads to
\[\lb{116}
 K_j(z)=\sum_{n=-\iy}^{+\iy}L_{nj}z^{\frac{2\pi\mathbf{i}n+\ln q_j}{\ln E}}.
\]

\subsection{Series for the density $p(x)$.}

Substituting \er{116} into \er{110}, we get
\[\lb{117}
 \Pi(z)=\sum_{j=1}^{+\iy}b_j\sum_{n=-\iy}^{+\iy}L_{nj}z^{\frac{2\pi\mathbf{i}n+\ln q_j}{\ln E}}.
\]
Then, by \er{101} and \er{117}, we have
\[\lb{118}
 p(x)=\sum_{j=1}^{+\iy}b_j\sum_{n=-\iy}^{+\iy}\frac{L_{nj}}{2\pi}\int_{-\iy}^{+\iy}(\ve+\mathbf{i}y)^{\frac{2\pi\mathbf{i}n+\ln q_j}{\ln E}}e^{(\ve+\mathbf{i}y)x}dy=\sum_{j=1}^{+\iy}b_j\sum_{n=-\iy}^{+\iy}\frac{L_{nj}x^{-1-^{\frac{2\pi\mathbf{i}n+\ln q_j}{\ln E}}}}{\G(-{\frac{2\pi\mathbf{i}n+\ln q_j}{\ln E}})}
\]
or
\[\lb{118a}
 p(x)=\sum_{j=1}^{+\iy}b_jx^{-1-\log_E q_j}A_j(-\log_Ex)
\]
with one-periodic functions
\[\lb{118b}
 A_j(x):=\sum_{n=-\iy}^{+\iy}\frac{L_{nj}e^{2\pi\mathbf{i}nx}}{\G(-{\frac{2\pi\mathbf{i}n+\ln q_j}{\ln E}})}.
\]
Identity \er{118a} coincides with \er{p_main}. In the last equality of \er{118}, we use the identity
\[\lb{118c}
 \frac1{\G(z)}=\frac1{2\pi}\int_{-\iy}^{+\iy}e^{\ve+\mathbf{i}y}(\ve+\mathbf{i}y)^{-z}dy,
\]
mentioned in, e.g., \cite{WW} p. 254 Example 12.2.6.

\subsection{Series for $\vp_{nj}$.} From \er{111}, we have
\[\lb{119}
 \Phi_j(1-z)\sim K_j(z),\ \ \ z\to0,
\]
at least for $\Re z>0$. This is true because
\[\lb{119a}
 P_{r_t}\circ...\circ P_{r_1}\lt(1-\frac{z}{E^t}\rt)\sim 1-z,\ \ \ z\to0,
\]
since $P_{r}(1)=1$ and $P_{r}'(1)=E$ for all $r$. Thus, due to \er{116}, we obtain
\[\lb{120}
 \Phi_j(z)\sim\sum_{n=-\iy}^{+\iy}L_{nj}(1-z)^{\frac{2\pi\mathbf{i}n+\ln q_j}{\ln E}}=\sum_{n=-\iy}^{+\iy}L_{nj}\sum_{m=0}^{+\iy}(-1)^m\binom{\frac{2\pi\mathbf{i}n+\ln q_j}{\ln E}}{m}z^m,\ \ \ z\to1,
\]
at least for $\Re z<1$. Using the asymptotic of the binomial coefficients, see, e.g., \cite{TE} and \cite{FO},
\[\lb{121}
 \binom{\a}{m}\sim\frac{(-1)^m}{\G(-\a)m^{\a+1}}S_{\a}(m),\ \ \ S_{\a}(m):=1+\frac{s_1(\a)}{m}+\frac{s_2(\a)}{m^2}+...,\ \ \ m\to+\iy,
\]
with
\[\lb{121a}
 s_1(\a):=\frac{\a(\a+1)}{2},\ \ \ s_2(\a):=\frac{\a(\a+1)(\a+2)(3\a+1)}{24},\ \ \ ...,
\]
and comparing \er{120} with \er{103}, we obtain
\[\lb{122}
 \vp_{mj}\sim m^{-1-\log_Eq_j}\lt(A_j(-\log_Em)+\frac{B_{1j}(-\log_Em)}{m}+\frac{B_{2j}(-\log_Em)}{m^2}+...\rt),
\]
where the corresponding one-periodic functions are defined by \er{118b} and by
\[\lb{123}
 B_{kj}(x):=\sum_{n=-\iy}^{+\iy}\frac{s_k(\frac{2\pi\mathbf{i}n+\ln q_j}{\ln E})L_{nj}e^{2\pi\mathbf{i}nx}}{\G(-{\frac{2\pi\mathbf{i}n+\ln q_j}{\ln E}})}.
\]
Identity \er{vphi_main} follows from \er{122}.

{\section{Examples}\lb{sec2}}

For most examples, Embarcadero Delphi Rad Studio Community Edition and the library NesLib.Multiprecision is used. This software provides a convenient environment for programming and well-functioning basic functions for high-precision computations. All the algorithms related to the article's subject, including efficient parallelization, are developed by the author (AK). We use 128-bit precision instead of the standard 64-bit double precision.

We consider the cases with only two polynomials
\[\lb{200}
 P_0(z)=pz+(1-p)z^2,\ \ \ P_1(z)=qz+(1-q)z^3,
\]
representing two different probability-generating functions. The expectation should be the same
\[\lb{201}
 E=p+2(1-p)=q+3(1-q),
\]
which gives
\[\lb{202}
 q=\frac{1+p}2.
\]
Thus, we have only one free parameter $p\in(0,1)$. We assume that $P_0$ and $P_1$ appear in the random environment with the same probability $1/2$. The functions
\[\lb{203}
 \Pi_t(z)=\sum_{{\bf r}\in\{0,1\}^t}2^{-t}P_{r_1}\circ...\circ P_{r_t}\lt(1-\frac{z}{E^t}\rt)
\]
tend to $\Pi(z)$, see \er{100}, when $t\to\iy$. In computations, we set $t=19$, which is already quite expensive. For the integral approximation of the density, see \er{101}, we use the trapezoidal rule for
\[\lb{204}
 \tilde p(x)=\Re\lt(\frac1{\pi}\int_{0}^{500}\Pi_{19}(\mathbf{i}y)e^{\mathbf{i}yx}dy\rt)
\]
with the discretization step $1/100$.

The computation of $\vp_{nj}$ is based on \er{105}, which for two polynomials \er{200} becomes
\[\lb{205}
 \vp_{nj}=(p^j+q^j-p^n-q^n)^{-1}\lt(\sum_{i=1}^{2i\le n}B_2(i)\vp_{n-i,j}+\sum_{i=1}^{3i\le n}B_3(i)\vp_{n-2i,j}\rt)
\]
with
\[\lb{206}
 B_2(i)=\frac{(n-i)!}{(n-2i)!\cdot i!}p^{n-2i}(1-p)^{i},\ \ \ B_3(i)=\frac{(n-2i)!}{(n-3i)!\cdot i!}q^{n-3i}(1-q)^{i}.
\]
Taking the logarithm of both identities in \er{206} and then the exponent improves the accuracy of computations. Also, having the precomputed values of $\ln i!$ for $i=1,...,I$ with some large $I$ speeds up computations of \er{205} essentially. Then, we use \er{108} along with \er{205} to calculate the approximation \er{approx2}.

\begin{figure}
	\centering
	\begin{subfigure}[b]{0.75\textwidth}
		\includegraphics[width=\textwidth]{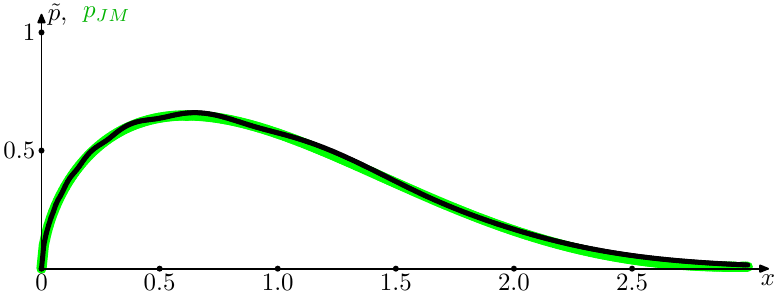}
		\caption{$p=0.2$}
		\label{fig1a}
	\end{subfigure}
	\hfill
	\begin{subfigure}[b]{0.75\textwidth}
		\includegraphics[width=\textwidth]{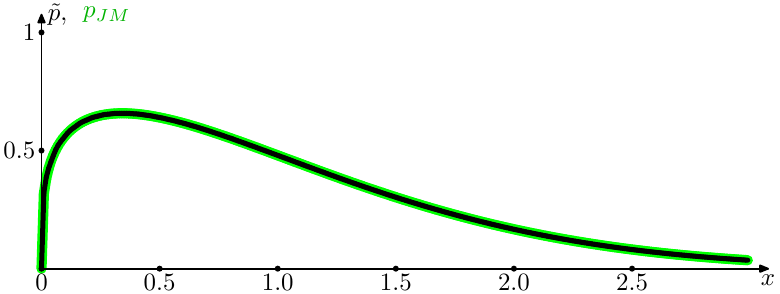}
		\caption{$p=0.4$}
		\label{fig1b}
	\end{subfigure}
	\hfill
	\begin{subfigure}[b]{0.75\textwidth}
		\includegraphics[width=\textwidth]{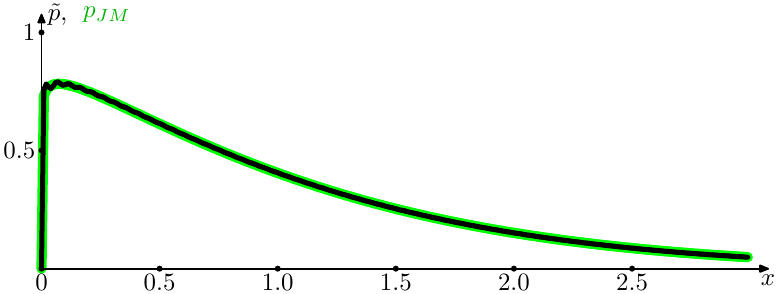}
		\caption{$p=0.6$}
		\label{fig1c}
	\end{subfigure}
	\caption{For the case \er{200}, the approximations of the density \er{101} are computed by \er{204} and \er{approx2} with $M=1000$ and $J=19$.}\label{fig1}
\end{figure}

The results of computations are presented in Fig. \ref{fig1}. The calculation of \er{approx2} is much faster than \er{204}. The first calculation is almost real-time, while the second may take up to 20 minutes on a Core i9-13900h. The oscillations of black curves seem spurious because the real oscillations are smaller. The reason for the spurious oscillations is a limited integration interval - the same happens for the classical case of a constant environment. 

{\section{Conclusion}\lb{sec3}}

We have derived a complete left tail asymptotic for the density of {\it martingale limit} for Galton-Watson processes with constant expectations in random environments. Based on the asymptotic, we obtained an easily computed approximation for the density, while the density itself is difficult to calculate. The following questions remain open: What happens in the case of non-constant mathematical expectation? What does the asymptotic of the right tail look like?

\section*{Ethical Statements}

The author has no relevant financial or non-financial interests to disclose. Only standard support for general research work without special interests is listed in the Acknowledgements section below.

The author has no competing interests to declare that are relevant to the content of this article.

The author certifies that he has no affiliations with or involvement in any organization or entity with any financial interest or non-financial interest in the subject matter or materials discussed in this manuscript.

The author has no financial or proprietary interests in any material discussed in this article.

\section*{Acknowledgements} 
This paper is a contribution to the project M3 of the Collaborative Research Centre TRR 181 "Energy Transfer in Atmosphere and Ocean" funded by the Deutsche Forschungsgemeinschaft (DFG, German Research Foundation) - Projektnummer 274762653. 

\bibliographystyle{abbrv}

\end{document}